\documentclass[11pt]{amsart}

\usepackage[frame,cmtip,arrow,matrix,line,graph,curve]{xy}
\usepackage{amsmath,amssymb, amscd}
\usepackage{amsfonts,amsxtra,psfig,amsthm}
\usepackage{mathrsfs}
\usepackage{eucal}
\usepackage{latexsym}

\numberwithin{equation}{section}

\newtheorem{theorem}{Theorem}[section]
\newtheorem{proposition}[theorem]{Proposition}
\newtheorem{corollary}[theorem]{Corollary}
\newtheorem{lemma}[theorem]{Lemma}

\newtheorem{question}[theorem]{Question}
\newtheorem{remit}[theorem]{Remark}
\newtheorem{assumpt}[theorem]{Assumption}
\newtheorem{definit}[theorem]{Definition}

\newenvironment{definition}{\begin{definit}\rm}{\end{definit}}

\newenvironment{remark}{\begin{remit}\rm}{\end{remit}}

\newcommand{\pp}{\mathbb{P}}
\newcommand{\qq}{\mathbb{Q}}
\newcommand{\cc}{\mathbb{C}}

\newcommand{\zz}{\mathbb{Z}}

\newcommand{\Hom}{\mathrm{Hom}}

\newcommand{\Gr}{\mathrm{Gr}}
\newcommand{\Ext}{\mathrm{Ext}}

\newcommand{\codim}{\mathrm{codim \,}}

\newcommand{\git}{/\!\!/}

\newcommand{\cA}{\mathcal{A} }
\newcommand{\cB}{\mathcal{B} }

\newcommand{\cE}{\mathcal{E} }

\newcommand{\cM}{\mathcal{M} }

\newcommand{\cL}{\mathcal{L} }
\newcommand{\cO}{\mathcal{O} }

\newcommand{\bM}{\mathbf{M} }

\begin{document}

\title[Nonexistence of a crepant resolution of moduli spaces]
{On the existence of a crepant resolution of some moduli spaces of
sheaves on an abelian surface}
\date{}

\author{Jaeyoo Choy and Young-Hoon Kiem}
\address{Dept of Mathematics, Seoul National University,
Seoul 151-747, Korea} \email{donvosco@math.snu.ac.kr}
\email{kiem@math.snu.ac.kr}
\thanks{Young-Hoon Kiem was partially supported by a KOSEF grant
R01-2003-000-11634-0}


\begin{abstract}
Let $J$ be an abelian surface with a generic ample line bundle
$\cO_J(1)$. For $n\ge 1$, the moduli space $M_J(2,0,2n)$ of
$\cO(1)$-semistable sheaves $F$ of rank $2$ with Chern classes
$c_1(F)=0$, $c_2(F)=2n$ is a singular projective variety, endowed
with a holomorphic symplectic structure on the smooth locus. In
this paper, we show that there does not exist a crepant resolution
of $M_J(2,0,2n)$ for $n\ge 2$. This certainly implies that there
is no symplectic desingularization of $M_J(2,0,2n)$ for $n\ge 2$.
\end{abstract}
\maketitle

\section{Introduction}

An irreducible symplectic manifold $X$ is a compact simply
connected complex manifold, endowed with a nondegenerate
holomorphic 2-form $\omega$ which spans $H^0(\Omega^2_X)$. By the
Bogomolov decomposition \cite{Bea83}, irreducible symplectic
manifolds are building blocks of K\"ahler manifolds in the sense
that for any compact K\"ahler manifold with trivial first Chern
class there is an \'etale cover from the product of tori,
Calabi-Yau manifolds and irreducible symplectic manifolds, which
we call the Bogomolov factors. Two standard series of examples
were provided by Beauville: Hilbert schemes of points on K3
surfaces and the generalized Kummer varieties \cite{Bea83}.

Recently O'Grady proposed a strategy for finding new examples of
irreducible symplectic manifolds as follows \cite{ogrady,og00}:
\begin{enumerate}
\item Consider a \emph{singular} moduli space  $M(r,c_1,c_2)$ of \emph{semistable}
sheaves on a K3 or abelian surface $S$ of rank $r$ with Chern
classes $c_1, c_2\in H^*(S,\zz)$. By Mukai's theorem \cite{Muk84},
there is a symplectic form, called the Mukai form, on the open
subset of stable sheaves $M(r,c_1,c_2)^s$.
\item Find a desingularization $\widetilde{M}(r,c_1,c_2)$ of
$M(r,c_1,c_2)$ on which the Mukai
form extends everywhere without degeneration.
\item Look at the Bogomolov factors for a new
irreducible symplectic manifold.\end{enumerate} Actually, O'Grady
successfully implemented his program in two cases and found new
irreducible symplectic manifolds of (complex) dimension 10 and 6
respectively:
\begin{enumerate} \item a symplectic desingularization of the moduli space $M_{K3}(2,0,4)$ of
rank $2$ semistable sheaves on a K3 surface with Chern classes
$c_1=0, c_2=4$  (see \cite{ogrady}), \item a Bogomolov factor of a
symplectic desingularization of the moduli space $M_{Ab}(2,0,2)$
of rank $2$ semistable sheaves on an abelian surface with Chern
classes $c_1=0, c_2=2$ (see \cite{og00}).\end{enumerate} A natural
question raised by O'Grady asks whether one can do the same with
$M_{K3}(2,0,2m)$ with $m\ge 3$ or $M_{Ab}(2,0,2n)$ with $n\ge 2$,
i.e.
\begin{question}\label{quest1} Does there exist a symplectic desingularization
of $M_{K3}(2,0,2m)$ with $m\ge 3$ or $M_{Ab}(2,0,2n)$ with $n\ge
2$?\end{question}
 In \cite{CK04,KL04}, it was proved that unfortunately the answer
is NO for the K3 case: there is no symplectic desingularization of
$M_{K3}(2,0,2m)$ for $m\ge 3$. However, the question remains open
for $M_{Ab}(2,0,2n)$ with $n\ge 2$. The purpose of this paper is
to show that the answer to the above question is also NO for the
abelian case, i.e. there is no symplectic desingularization of
$M_{Ab}(2,0,2n)$ with $n\ge 2$.

Fix any integer $n\ge 2$. Let $J$ be a complex projective abelian
surface equipped with a generic ample divisor $\Theta$, which
satisfies
\begin{assumpt}\cite[(1.3)]{og00}
\label{assumption on J}  There is no divisor $A$ orthogonal to
$\Theta$ with $-2n\le A^2<0$.
\end{assumpt} This condition is satisfied if for instance the
N\'eron-Severi group is $NS(J)=\zz c_1(\Theta).$ Let
$$\bM=\bM_{2n}=M_J(2,0,2n)$$
denote the moduli space of $\Theta$-semistable sheaves $F$ on $J$
of rank 2 with $c_1(F)=0$, $c_2(F)=2n$ in $H^*(J,\zz)$. This is an
irreducible normal projective variety of dimension $8n+2$
(\cite{Gi77, Muk84,Yo01}) with Gorenstein
singularities.\footnote{A normal variety $X$ is Gorenstein if the
canonical divisor $K_X$ is a Cartier divisor \cite{Bat98}. In our
case, the stable part $\bM^s$ of $\bM$ is equipped with a
symplectic form \cite{Muk84} and $\codim (\bM-\bM^s)\ge 2$ (see
\S3). Therefore, the canonical divisor $K_\bM$ is zero and thus
$K_\bM$ is Cartier.} The main result of this paper is the
following.
\begin{theorem} \label{thm:main result}
If $n \geq 2$, there is no crepant resolution of $\bM _{2n}$.
\end{theorem}
By Mukai's theorem \cite{Muk84}, the canonical line bundle of
$\bM_{2n}$ is trivial and hence any desingularization of
$\bM_{2n}$ equipped with a holomorphic symplectic form is a
crepant resolution. So we deduce from Theorem \ref{thm:main
result} the following corollary which answers Question
\ref{quest1}.
\begin{corollary} \label{cor:O'Grady's conjecture} If $n\geq2$,
there is no symplectic desingularization of $\bM_{2n}$.
\end{corollary}

As in \cite{CK04}, the idea of the proof of Theorem \ref{thm:main
result} is to use properties of the stringy E-function
\cite{Bat98}. Recall that $\bM$ is a normal irreducible variety
with Gorenstein singularities. Also we will see in section
\ref{sec: Kirwan's desingularization} that the singularities are
terminal. Hence, the stringy E-function $E_{st}({\bf M};u,v)$ of
$\bM$ is a well-defined rational function in formal variables
$u,v$.  By Kontsevich's Theorem (Theorem \ref{thm:Batyrev's
result}), if there is a crepant resolution $\widetilde \bM$ of
$\bM $, then the stringy E-function of $\bM$ is equal to the
Hodge-Deligne polynomial (E-polynomial) of $\widetilde \bM$. In
particular, we deduce that the stringy E-function $E_{st}({\bf
M};u,v)$ must be a polynomial. Therefore, Theorem \ref{thm:main
result} is a consequence of the following.
\begin{proposition}\label{prop:stringy E-function test} The stringy
E-function $E_{st}({\bf M} ;u,v)$ is not a polynomial for
$n\geq2$.
\end{proposition}

In \cite{ogrady, og97}, O'Grady studies Kirwan's desingularization
$\widehat{\bM}$ of $\bM$ which is obtained as the result of three
blow-ups. We use O'Grady's analysis of Kirwan's desingularization
of $\bM$ to prove Proposition \ref{prop:stringy E-function test}.
In section \ref{sec: properties of E-functions} we recall some
properties of stringy E-functions and we prove Proposition
\ref{prop:stringy E-function test} in section \ref{sec:
Nonexistence of Crepant Resolution}. In section \ref{sec: Kirwan's
desingularization} we analyze Kirwan's desingularization.

After completion of this paper, Kaledin, Lehn and Sorger proved
nonexistence of symplectic desingularization for arbitrary rank by
a different method. See \cite{KLS}. We are grateful to C. Sorger
for delightful conversations about singular symplectic moduli
spaces at the Korea Institute of Advanced Study during a workshop
on Vector Bundles on Algebraic Curves organized by S. Ramanan and
J.-M. Hwang in April 2005.


\section{Some properties of Poincar\'{e} polynomials,
Hodge-Deligne polynomials and stringy E-functions} \label{sec:
properties of E-functions}

In this section we collect some facts that we shall use later.

For a topological space $V$, the Poincar\'e polynomial of $V$ is
defined by
\begin{equation} \label{eqn:Poincare polynomial}
P(V;z)=\sum_{i}(-1)^ib_i(V)z^i
\end{equation} where $b_i(V)$ is the $i$-th Betti number of $V$.

Next we recall the definition and basic facts about stringy
E-functions from \cite{Bat98,DL99}. Let $W$ be a normal
irreducible variety with at worst log-terminal singularities, i.e.
\begin{enumerate} \item W is $\qq$-Gorenstein;
\item for a resolution of singularities $\rho: V\to W$ such that
the exceptional locus of $\rho$ is a divisor $D$ whose irreducible
components $D_1,\cdots,D_r$ are smooth divisors with only normal
crossings, we have \[K_V=\rho^*K_W+\sum^r_{i=1} a_iD_i \] with
$a_i>-1$ for all $i$, where $D_i$ runs over all irreducible
components of $D$. The divisor $\sum^r_{i=1}a_iD_i$ is called the
\textit{discrepancy divisor}. \end{enumerate}

\begin{definition}\label{def2.1} For each subset $J\subset I=\{1,2,\cdots,r\}$, define
$D_J=\cap_{j\in J}D_j$, $D_\emptyset =V$ and $D^0_J=D_J-\cup_{i\in
I-J}D_i$. Then the stringy E-function of $W$ is defined by
\begin{equation} \label{eqn:stringy E-function}
E_{st}(W;u,v)=\sum_{J\subset I}E(D^0_J;u,v)\prod_{j\in
J}\frac{uv-1}{(uv)^{a_j+1}-1} \end{equation} where \[ E(Z;u,v) =
\sum_{p,q}\sum_{k\geq 0} (-1)^kh^{p,q}(H^k_c(Z;\cc))u^pv^q \] is
the Hodge-Deligne polynomial (= E-polynomial) for a variety
$Z$.\end{definition} Note that the Hodge-Deligne polynomials have
\begin{enumerate}
\item the additive property: $E(Z;u,v)=E(U;u,v)+E(Z-U;u,v)$ if $U$
is an open subvariety of $Z$; \item the multiplicative property:
$E(Z;u,v)=E(B;u,v)E(F;u,v)$ if $Z$ is a Zariski locally trivial
$F$-bundle over $B$. \end{enumerate}

By \cite[Theorem 6.27]{Bat98}, the function $E_{st}$ is
independent of the choice of a resolution (Theorem 3.4 in
\cite{Bat98}) and the following holds.
\begin{theorem} \label{thm:Batyrev's result} \cite[Theorem
3.12]{Bat98} Suppose $W$ is a $\qq$-Gorenstein algebraic variety
with at worst log-terminal singularities. If $\rho:V\to W$ is a
crepant desingularization (i.e. $\rho^*K_W=K_V$) then
$E_{st}(W;u,v)=E(V;u,v)$. In particular, $E_{st}(W;u,v)$ is a
polynomial.
\end{theorem}



\section{Kirwan's desingularization of $\bM$}
\label{sec: Kirwan's desingularization} In this section, we
analyze Kirwan's desingularization
 $$\rho:\widehat{\bM}\to\bM$$ constructed in
\cite[\S2.1]{og00}. Only Propositions \ref{prop:analysis on exc}
and \ref{discrepancy divisor} will be used in section \ref{sec:
Nonexistence of Crepant Resolution}.

 Let $C$ be a smooth irreducible projective
curve of genus 2 and $J=\mathrm{Pic}^0(C)$. Fix a Weierstrass
point $p_0$ of $C$ and let $\Theta$ be the image of the
Abel-Jacobi map $C\to J$ defined by $p\mapsto p-p_0$. In this
paper, we always suppose Assumption \ref{assumption on J} is
satisfied as in \cite[(1.3)]{og00}. This is obviously satisfied if
for instance the N\'{e}ron-Severi group is $NS(J) =\zz
c_1(\Theta)$. From now on, (semi)stability of a torsion-free sheaf
on $J$ means (semi)stability with respect to the ample divisor
$\Theta=:\cO_J(1)$. Let $J^{[n]}$ denote the Hilbert scheme of $n$
points in $J$ and $\hat{J}=\mathrm{Pic}^0(J)$.

Assumption \ref{assumption on J} is necessary for the following
(\cite[Lemma 2.1.2]{og00}):
\begin{lemma} A torsion-free sheaf $F$ of rank 2 with $c_1(F)=0$
on $J$ is strictly semistable if and only if $F$ fits into a short
exact sequence
$$0\to I_{Z_1}\otimes \xi_1\to F\to I_{Z_2}\otimes \xi_2\to 0$$
where $I_{Z_i}$ is the ideal sheaf of a zero dimensional subscheme
$Z_i\in J^{[n]}$ of length $n$ and $\xi_i\in
\mathrm{Pic}^0(J)=\hat{J}$.\end{lemma}

Now consider Simpson's construction of the moduli space
$\bM=\bM_{2n}$ \cite[\S1.1]{og97}. Let $Q$ be the closure of the
set of semistable points $Q^{ss}$ in the Quot-scheme whose
quotient by the natural $PGL(N)$ action is $\bM$ for some even
integer $N$. Then $Q^{ss}$ parameterizes semistable sheaves $F$
together with surjective homomorphisms $h:\cO^{\oplus N}\to F(k)$
which induces an isomorphism $\cc^N\cong H^0(F(k))$. Let
$\Omega_Q$ denote the subset of $Q^{ss}$ which parameterizes
sheaves of the form $(I_Z\otimes\xi)^{\oplus 2}$ for some $Z\in
J^{[n]}$ and $\xi\in \mathrm{Pic}^0(J)=\hat{J}$. Then $\Omega_Q$
is precisely the locus of closed orbits with maximal dimensional
stabilizers, isomorphic to $PGL(2)$ and the quotient of $\Omega_Q$
by $PGL(N)$ is $$\Omega:=\Omega_Q\git PGL(N)\cong J^{[n]}\times
\hat{J}.$$ Let $\Sigma_Q$ be the subscheme
 of $Q^{ss}$ which parameterizes sheaves of the form $(I_{Z_1}\otimes
 \xi_1)\oplus (I_{Z_2}\otimes
 \xi_2)$ for some $Z_1,Z_2\in J^{[n]}$ and $\xi_1,\xi_2\in
 \mathrm{Pic}^0(J)=\hat{J}$. Then $\Sigma_Q-\Omega_Q$ is precisely the
 locus of closed orbits with 1-dimensional stabilizers isomorphic
 to $\cc^*$. The quotient of $\Sigma_Q$ by $PGL(N)$ is
 $$\Sigma:=\Sigma_Q\git PGL(N)\cong (J^{[n]}\times \hat{J})\times
  (J^{[n]}\times \hat{J})/\zz_2$$
 where the $\zz_2$-action is the involution which interchanges the
 two components while $\Omega$ sits in $\Sigma$ as the diagonal.
 So we have a stratification of $\bM$:
$$\bM=\bM^s\sqcup (\Sigma-\Omega)\sqcup \Omega$$
where $\bM^s$ is the locus of stable sheaves which is smooth by
\cite{Muk84}. To obtain a desingularization of $\bM$ we blow up
$\bM$ along $\Omega$ and then along the proper transform of
$\Sigma$. The result of these two blow-ups is an orbifold and by
blowing up once more along the singular locus we get a smooth
model of $\bM$, which we call Kirwan's desingularization
\cite[\S2.1]{og00}, \cite[Proposition 1.8.3]{ogrady}.

For a detailed analysis of $\Omega_Q$ and $\Sigma_Q$, we need to
make some observations. To begin with, note that at each
$(Z,\xi)\in J^{[n]}\times \hat{J}$ the tangent space
$T_{J^{[n]}\times \hat J,(Z,\xi)}$ of $J^{[n]}\times \hat J\simeq
{M}_J(1,0,n)$ is canonically isomorphic to $\Ext^1(I_Z,I_Z)$ where
$I_Z$ is the ideal sheaf of $Z$. By the Yoneda pairing map and the
Serre duality, we have a skew-symmetric pairing
$$\omega:\Ext^1(I_Z,I_Z)\otimes \Ext^1(I_Z,I_Z) \to
\Ext^2(I_Z,I_Z)\cong \cc$$ which gives us a symplectic form
$\omega$ on the tangent bundle $T_{J^{[n]}\times \hat J}$ by
\cite[Theorem 0.1]{Muk84}.

Let $W=sl(2)^\vee\cong sl(2)\cong \cc^3$. The adjoint action of
$PGL(2)$ on $W$ gives us an identification $SO(W)\cong PGL(2)$
(\cite[\S1.5]{og97}). For a symplectic vector space $(V,\omega)$,
let $\Hom^\omega(W,V)$ be the space of homomorphisms from $W$ to
$V$ whose image is isotropic, i.e. the restriction of $\omega$ to
the image is trivial. Let $\Hom^\omega(W,T_{J^{[n]}\times \hat J})
$ be the bundle over $J^{[n]}\times \hat J$ whose fiber over
$(Z,\xi)\in J^{[n]}\times \hat J$ is
$\Hom^\omega(W,T_{J^{[n]}\times \hat J,(Z,\xi)})$. As an algebraic
vector bundle, $T_{J^{[n]}\times \hat J}$ is a Zariski locally
trivial bundle. By elementary linear algebra, we can furthermore
find local trivializations so that the symplectic form $\omega$ is
given by a constant skew-symmetric matrix on each open set.
Therefore, the bundle $\Hom^\omega(W,T_{J^{[n]}\times \hat J})$ is
Zariski locally trivial. Let $\Hom_k^\omega(W,T_{J^{[n]}\times
\hat J})$ be the subbundle of $\Hom^\omega(W,T_{J^{[n]}\times \hat
J})$ of rank $\leq k$ elements in $\Hom^\omega(W,T_{J^{[n]}\times
\hat J})$. Also let $\Gr^\omega(3,T_{J^{[n]}\times \hat J})$ be
the relative Grassmannian of isotropic 3-dimensional subspaces in
$T_{J^{[n]}\times \hat J}$ and let $\cB$ denote the tautological
rank 3 bundle on $\Gr^\omega(3,T_{J^{[n]}\times \hat J})$.
Obviously these bundles are all Zariski locally trivial as well.

Let $\pp\Hom^\omega(W,T_{J^{[n]}\times \hat J})$ (resp.
$\pp\Hom_k^\omega(W,T_{J^{[n]}\times \hat J})$) be the
projectivization of $\Hom^\omega(W,T_{J^{[n]}\times \hat J})$
(resp. $\Hom_k^\omega(W,T_{J^{[n]}\times \hat J})$). Likewise, let
$\pp\Hom(W,\cB)$ and $\pp\Hom_k(W,\cB)$ denote the
projectivizations of the bundles $\Hom(W,\cB)$ and
$\Hom_k(W,\cB)$. Note that there are obvious forgetful maps
\begin{eqnarray*}f:\pp\Hom(W,\cB)\to\pp\Hom^\omega(W,T_{J^{[n]}\times \hat J})\
 \mbox{\rm and}\\
 f_k:\pp\Hom_k(W,\cB)\to\pp\Hom_k^\omega(W,T_{J^{[n]}\times \hat J})
 \end{eqnarray*}
Since the pull-back of the defining ideal of
$\pp\Hom_1^\omega(W,T_{J^{[n]}\times \hat J})$ is the ideal of
$\pp\Hom_1(W,\cB)$ (both are actually given by the determinants of
$2\times 2$ minor matrices), $f$ gives rise to a map between
blow-ups
$$\overline{f}:Bl_{\pp\Hom_1(W,\cB)}\pp\Hom(W,\cB)\to
Bl_{\pp\Hom_1^\omega(W,T_{J^{[n]}\times \hat
J})}\pp\Hom^\omega(W,T_{J^{[n]}\times \hat J}).$$ Let us denote
$Bl_{\pp\Hom_1(W,\cB)}\pp\Hom(W,\cB)$ by $Bl^\cB$ and
$Bl_{\pp\Hom_1^\omega(W,T_{J^{[n]}\times \hat
J})}\pp\Hom^\omega(W,T_{J^{[n]}\times \hat J})$ by $Bl^T$. We
denote the proper transform of $\pp\Hom_2(W,\cB)$ in $Bl^\cB$ by
$Bl_2^\cB$ and the proper transform of
$\pp\Hom_2^\omega(W,T_{J^{[n]}\times \hat J})$ by $Bl_2^T$. Since
$Bl_2^\cB$ is a smooth divisor which is mapped onto $Bl_2^T$ and
the pull-back of the defining ideal of $Bl_2^T$ is the ideal sheaf
of $Bl_2^\cB$, $\overline{f}$ lifts to
\begin{equation}\label{eq4.-2}
\hat{f}:Bl^\cB \to Bl_{Bl_2^T}Bl^T.\end{equation} By \cite[\S3.1
IV]{og97}, $\hat f$ is an isomorphism on each fiber over
$J^{[n]}$, so in particular $\hat f$ is bijective. Therefore,
$\hat f$ is an isomorphism.

Note that $\pp\Hom(W,\cB)\git SO(W)$ (resp. $\pp\Hom_k(W,\cB)\git
SO(W)$) is isomorphic to the space of conics $\pp(S^2\cB)$ (resp.
rank $\leq k$ conics $\pp(S^2_k\cB)$)  where the $SO(W)$-quotient
map is given by $\alpha\mapsto\alpha\circ\alpha^t$
(\cite[\S3.1]{og97}). Let $\hat \pp(S^2\cB)\cong
Bl_{\pp(S^2_1\cB)}\pp(S^2\cB)$ denote the blow-up along the locus
of rank 1 conics. Then $Bl^\cB\git SO(W)$ is canonically
isomorphic to $\hat \pp(S^2\cB)$ by \cite[Lemma 3.11]{k2}. Since
$\cB$ is Zariski locally trivial, so is $\hat \pp(S^2\cB)$ over
$\Gr^\omega(3,T_{J^{[n]}\times \hat J})$.

Now we can give a more precise description of $\Omega_Q$ as
follows. Let $\cL$ be a universal rank 1 sheaf over
$(J^{[n]}\times \hat J)\times   J=M_J(1,0,n)\times J$ such that
$\cL|_{(Z,\xi)\times J}$ is isomorphic to $I_Z\otimes \xi$. By
\cite[Theorem 10.2.1]{HL97},  the tangent bundle $T_{J^{[n]}\times
\hat J}$ is in fact isomorphic to $\cE xt^1_{J^{[n]}\times \hat J
}(\cL,\cL)$. Let $p:(J^{[n]}\times \hat J)\times J\to
J^{[n]}\times \hat J$ be the projection onto $J^{[n]}\times \hat J
$ and $p_J$ be the projection onto $J$. By tensoring with the
pull-back of $\cO_J(k)$ for suitable $k$, $p_*\cL(k)$ is a vector
bundle of rank $N/2$ where $\cL(k):=\cL\otimes p_J^*\cO_J(k)$. Let
\begin{equation}\label{eq4.-1}
 q:\pp \mathrm{Isom}(\cc^N, p_*\cL(k)\oplus
 p_*\cL(k))\to J^{[n]}\times \hat J
 \end{equation}
be the $PGL(N)$-bundle over $J^{[n]}\times \hat J$ whose fiber
over $(Z,\xi)$ is $$\pp \mathrm{Isom}(\cc^N, H^0((I_Z\otimes
\xi)^{\oplus 2}\otimes\cO_J(k))).$$ Note that the standard action
of $PGL(N)$ on $\cc^N$ commutes with the obvious action of
$PGL(2)\cong SO(W)$ on $p_*\cL(k)\oplus p_*\cL(k)$.
\begin{lemma}\label{4.1} (1) $\Omega_Q\cong \pp
\mathrm{Isom}(\cc^N, p_*\cL(k)\oplus p_*\cL(k))\git SO(W).$\\
(2) Via the above isomorphism, the normal cone of $\Omega_Q$ in
$Q^{ss}$ is $$q^*\mathrm{Hom}^{\omega}(W,T_{J^{[n]}\times \hat
J})\git SO(W)\to \pp \mathrm{Isom}(\cc^N, p_*\cL(k)\oplus
p_*\cL(k))\git SO(W)$$ whose fiber is
$\mathrm{Hom}^{\omega}(W,T_{J^{[n]} \times \hat J,(Z,\xi)})$.
\end{lemma}
\begin{proof}
(1) This is standard and we omit the proof.

(2) Let $\cO^{\oplus N} \twoheadrightarrow \cE$  denote the
universal quotient sheaf on $Q^{ss}\times J$. The Kodaira-Spencer
map associated to $\cE$ restricted to $\Omega_Q$ gives us a map
from the tangent sheaf $T_{Q^{ss}} |_{\Omega_Q}$ to the sheaf $\cE
xt^1_{\Omega_Q}(\cE,\cE)$ whose kernel is the tangent sheaf of the
orbits. Via the isomorphism of (1), we have a map
$$\delta:\pp
\mathrm{Isom}(\cc^N, p_*\cL(k)\oplus p_*\cL(k))\to \pp
\mathrm{Isom}(\cc^N, p_*\cL(k)\oplus p_*\cL(k))\git SO(W)\cong
\Omega_Q.$$ From the proof of (1) above, the pull-back of $\cE$ by
 $\delta$ is isomorphic to $(q\times 1)^*(\cL(k)\oplus
\cL(k))\otimes H$ and thus the vector bundle $\delta^*\cE
xt^1_{\Omega_Q}(\cE,\cE)$ is isomorphic to
$$q^*\cE xt^1_{J^{[n]}\times \hat J}(\cL, \cL)\otimes gl(2)\cong
q^*T_{J^{[n]}\times \hat J}\otimes gl(2).$$ The pull-back of the
tangent sheaf of $J^{[n]}\times \hat J$ sits in it as
$q^*T_{J^{[n]}\times \hat J}\otimes \left(\begin{matrix}1&0\\
0&1\end{matrix}\right)$ and thus the pull-back by $\delta$ of the
normal sheaf to $\Omega_Q$ is isomorphic to
$$q^*T_{J^{[n]}\times \hat J}\otimes sl(2)\cong q^*\mathrm{Hom}
(W,T_{J^{[n]}\times \hat J}).$$ By \cite{og97} (1.4.10), the
normal cone is the same as the Hessian cone fiberwisely. Since the
normal cone is contained in the Hessian cone, the normal cone is
equal to the Hessian cone which is the inverse image of zero by
the Yoneda square map $\Upsilon:\cE xt^1_{\Omega_Q}(\cE,\cE)\to
\cE xt^2_{\Omega_Q}(\cE,\cE)$. It is elementary to see that
$\delta^*\Upsilon^{-1}(0)$ is precisely
$q^*\mathrm{Hom}^\omega(W,T_{J^{[n]}\times \hat J}).$ Since
$SO(W)$ acts freely we obtain (2).
\end{proof}

 Let
$\pi_R:R\to Q^{ss}$ be the blow-up of $Q^{ss}$ along $\Omega_Q$.
Let $\Omega_R$ be the exceptional divisor of $\pi_R$ and
$\Sigma_R$ be the proper transform of $\Sigma_Q$. By the above
lemma, we have
\begin{equation}\label{eq4.0}\Omega_R\cong
 q^*\pp\mathrm{Hom}^{\omega}(W,T_{J^{[n]}\times \hat J})\git
 SO(W).
 \end{equation}
The following lemma is an easy exercise.
 \begin{lemma}\label{4.2} (1) The locus of points
in $\pp\mathrm{Hom}^\omega(W,T_{J^{[n]}\times \hat
J,(Z,\xi)})^{ss}$ whose stabilizer is 1-dimensional by the action
of $SO(W)$ is precisely $$\pp \mathrm{Hom}^\omega
_1(W,T_{J^{[n]}\times \hat J,(Z,\xi)})^{ss}.$$ (2) The locus of
nontrivial stabilizers is $\pp \mathrm{Hom}^\omega
_2(W,T_{J^{[n]}\times \hat J,(Z,\xi)})^{ss}$ and the stabilizers
are isomorphic to $\zz_2$ or $\cc^*$.\end{lemma} Let
\begin{equation}\label{eq4.1}
 \Delta_R=q^*\pp\mathrm{Hom}^{\omega}_2(W,T_{J^{[n]}\times \hat J})\git
 SO(W).
 \end{equation}
Note that $\Sigma_Q-\Omega_Q$ is precisely the locus of points in
$Q^{ss}$ whose stabilizer is isomorphic to $\cc^*$ and hence
$\Sigma_R^{ss}$ is precisely the locus of points in $R^{ss}$ with
1-dimensional stabilizers by \cite{k2}. Therefore we have the
following from Lemma \ref{4.2}.
\begin{corollary}\label{4.3}
$\Sigma_R^{ss}\cap
\Omega_R=q^*\pp\mathrm{Hom}^{\omega}_1(W,T_{J^{[n]}\times \hat
J})^{ss}\git SO(W).$\end{corollary}

 We have an explicit description of
$\Sigma_R^{ss}$ which is parallel to \cite[\S1.7 III]{og97} as
follows. Let
$$\beta:\mathcal{J}^{[n]}\to (J^{[n]}\times \hat J)\times (J^{[n]}\times \hat J)$$
be the blow-up along the diagonal and let $\mathcal{J}^{[n]}_0=
(J^{[n]}\times \hat J)\times (J^{[n]}\times \hat J)-
\mathbf{\Delta}$ where $\mathbf{\Delta}$ is the diagonal. Let
$\cL_1$ (resp. $\cL_2$) be the pull-back to
$\mathcal{J}^{[n]}\times J$ of the universal sheaf $\cL\to
(J^{[n]}\times \hat J)\times J$ by  $p_{13}\circ (\beta\times 1)$
(resp. $p_{23}\circ (\beta\times 1)$) where $p_{ij}$ is the
projection onto the first (resp. second) and third components. Let
$p:\mathcal{J}^{[n]}\times J\to \mathcal{J}^{[n]}$ be the
projection onto the first component. Then $p_*\cL_1(k)\oplus
p_*\cL_2(k)$ is a vector bundle of rank $N$. Let
$$q:\pp\mathrm{Isom}(\cc^N,p_*\cL_1(k)\oplus p_*\cL_2(k))\to \mathcal{J}^{[n]}$$
be the $PGL(N)$-bundle. There is an action of $O(2)$ on
$\pp\mathrm{Isom}(\cc^N,p_*\cL_1(k)\oplus p_*\cL_2(k))$. The
following lemma is obtained by (a proof parallel to) \cite{og97}
(1.7.10) and (1.7.1).
\begin{lemma}\label{4.4}
(1) $\Sigma_R^{ss}\cong\pp\mathrm{Isom}(\cc^N,p_*\cL_1(k)\oplus
p_*\cL_2(k))\git O(2)$\\
(2) The normal cone of $\Sigma_R^{ss}$ in $R^{ss}$ is a locally
trivial bundle over $\Sigma_R^{ss}$ with fiber the cone over a
smooth quadric in $\pp^{4n-1}$.
\end{lemma}
In fact we can give a more explicit description of the normal cone
when restricted to $\Sigma_R^0:=\Sigma_R^{ss}-\Omega_R$. Similarly
as in the proof of Lemma \ref{4.1}, the normal sheaf to
$\Sigma_R^0$ is isomorphic to the vector bundle (of rank $4n$)
\begin{equation}\label{eqn:4.2}
 q^*[\cE xt^1_{\mathcal{J}^{[n]}_0}(\cL_1,\cL_2)\oplus \cE
 xt^1_{\mathcal{J}^{[n]}_0}(\cL_2,\cL_1)]\git O(2)
 \end{equation}
over $\pp\mathrm{Isom}(\cc^N,p_*\cL_1(k)\oplus p_*\cL_2(k))\git
O(2)$ where $O(2)$ acts as follows: if we realize $O(2)$ as the
subgroup of $PGL(2)$
generated by $$SO(2)=\{\theta_\alpha=\left(\begin{matrix}\alpha&0\\
0&\alpha^{-1}\end{matrix}\right)\}/\{\pm Id\},\qquad
\tau=\left(\begin{matrix} 0&1\\1&0\end{matrix}\right)$$
$\theta_\alpha$ multiplies $\alpha$ (resp. $\alpha^{-1}$) to
$\cL_1$ (resp. $\cL_2$) and $\tau$ interchanges $\cL_1$ and
$\cL_2$. The normal cone is the inverse image
$q^*\Upsilon^{-1}(0)$ of zero in terms of the Yoneda pairing
\begin{equation}\label{eq4.3}\Upsilon:\cE
xt^1_{\mathcal{J}^{[n]}_0}(\cL_1,\cL_2)\oplus \cE
xt^1_{\mathcal{J}^{[n]}_0}(\cL_2,\cL_1)\to \cE
xt^2_{\mathcal{J}^{[n]}_0}(\cL_1,\cL_1).\end{equation}

Let $\pi_S:S\to R^{ss}$ be the blow-up of $R^{ss}$ along
$\Sigma_R^{ss}$. Let $\Sigma_S$ be the exceptional divisor of
$\pi_S$ and $\Omega_S$ (resp. $\Delta_S$) be the proper transform
of $\Omega_R$ (resp. $\Delta_R$). By \eqref{eq4.3}, we have
\begin{equation}\label{eq4.4}
 {\Sigma_S}
 |_{\pi_S^{-1}(\Sigma_R^0)}\cong q^*\pp \Upsilon^{-1}(0)\git
 O(2)\subset q^*\pp[\cE
 xt^1_{\mathcal{J}^{[n]}_0}(\cL_1,\cL_2)\oplus \cE
 xt^1_{\mathcal{J}^{[n]}_0}(\cL_2,\cL_1)]\git O(2).
 \end{equation}
By (a proof parallel to) \cite[(1.8.10)]{og97}, $S^s=S^{ss}$ and
$S^s$ is smooth. The quotient $S\git PGL(N)$ has only
$\zz_2$-quotient singularities along $\Delta_S\git PGL(N)$.

Finally let $\pi_T:T\to S^{s}$ be the blow-up of $S^s$ along
$\Delta_S^s$. Let $\Delta_T$ be the exceptional divisor of $\pi_T$
and $\Omega_T$ (resp. $\Sigma_T$) be the proper transform of
$\Omega_S$ (resp. $\Sigma_S$). Since $\Omega_T^{s} $,
$\Sigma_T^{s}$ and $\Delta_T^{s}$ are smooth divisors with finite
stabilizers $T\git PGL(N)$ is nonsingular and this is Kirwan's
desingularization $$\rho:\widehat{\bM}\to \bM.$$ The quotients
$\Omega_T\git PGL(N)$, $\Sigma_T\git PGL(N)$ and $\Delta_T\git
PGL(N)$ are denoted by $D_1=\hat\Omega$, $D_2=\hat\Sigma$ and
$D_3=\hat\Delta$ respectively.

We are ready to describe all the intersections of the smooth
divisors $D_1$, $D_2$ and $D_3$. Let $\hat{\pp}^5$ be the blow-up
of $\pp^5$ (projectivization of the space of $3\times 3$ symmetric
matrices) along $\pp^2$ (the locus of rank 1 matrices). For  a
symplectic vector space $(\cc^{2n},\omega)$, $\Gr^{\omega}(k,2n)$
denotes the Grassmannian of $k$-dimensional subspaces of
$\cc^{2n}$, isotropic with respect to the symplectic form $\omega$
(i.e. the restriction of $\omega$ to the subspace is zero).

\begin{proposition}\label{prop:analysis on exc}
Let $n\geq 2$.

(1) $D_1$ is a $\hat{\pp}^5$-bundle over a
$\Gr^\omega(3,2n+2)$-bundle over $J^{[n]}\times \hat J$.

(2) $D_2^0$ is a free $\zz_2$-quotient of a Zariski locally
trivial $I_{2n-1}$-bundle over $\mathcal{J}_0^{[n]}=(J^{[n]}\times
\hat J)\times (J^{[n]}\times \hat J)-\mathbf{\Delta}$ where
$\mathbf{\Delta}$ is the diagonal and $I_{2n-1}$ is the incidence
variety given by
\[ I_{2n-1}=\{(p,H)\in \pp^{2n-1}\times \breve{\pp}^{2n-1}| p\in
H\}.
\]

(3) $D_3$ is a $\pp^{2n-2}$-bundle over a Zariski locally trivial
$ \pp^2$-bundle over a Zariski locally trivial
$\Gr^\omega(2,2n+2)$-bundle over $J^{[n]}\times \hat J$.

(4) $D_1\cap D_2$ is a $\pp^2\times \pp^2$-bundle over a
$\Gr^\omega(3,2n+2)$-bundle over $J^{[n]}\times \hat J$.

(5) $D_2\cap D_3$ is a $\pp^{2n-2}$-bundle over a Zariski locally
trivial $ \pp^1$-bundle over a Zariski locally trivial
$\Gr^\omega(2,2n+2)$-bundle over $J^{[n]}\times \hat J$.

(6) $D_1\cap D_3$ is a $ \pp^2\times \pp^2$-bundle over a
$\Gr^\omega(3,2n+2)$-bundle over $J^{[n]}\times \hat J$.

(7) $D_1\cap  D_2 \cap  D_3$ is a $\pp^1\times \pp^2$-bundle over
a $\Gr^\omega(3,2n+2)$-bundle over $J^{[n]}\times \hat J$. \\
All the above bundles except in (2), (3) and (5)  are Zariski
locally trivial. Moreover, $D_i$ ($i=1,2,3$) are smooth divisors
such that $D_1\cup D_2\cup D_3$ is normal crossing.
\end{proposition}

\proof (1)  By \eqref{eq4.0} and Corollary \ref{4.3}, $\Omega_S$
is the blow-up of
$$q^*\pp\mathrm{Hom}^{\omega}(W,T_{J^{[n]}\times \hat J})\git SO(W)\text{ along
}q^*\pp\mathrm{Hom}^{\omega}_1(W,T_{J^{[n]}\times \hat J})\git
SO(W).$$ By \eqref{eq4.1}, $\Omega_T$ is the blow-up of $\Omega_S$
along the proper transform of
$$q^*\pp\mathrm{Hom}^{\omega}_2(W,T_{J^{[n]}\times \hat J})\git SO(W)$$ and
$D_1=\hat\Omega$ is the quotient of $\Omega_T$ by the action of
$PGL(N)$. Since the action of $PGL(N)$ commutes with the action of
$SO(W)$, $D_1$ is in fact the quotient by $SO(W)\times PGL(N)$ of
the variety obtained from
$q^*\pp\mathrm{Hom}^{\omega}(W,T_{J^{[n]}\times \hat J})$ by two
blow-ups. So $D_1$ is also the consequence of taking the quotient
by $PGL(N)$ first and then the quotient by $SO(W)$ second. Since
$q$ \eqref{eq4.-1} is a principal $PGL(N)$-bundle, the result of
the first quotient is just $Bl_{Bl_2^T}Bl^T$ in \eqref{eq4.-2}
which is isomorphic to $Bl^\cB$. If we take further the quotient
by $SO(W)$, then as discussed above the result is $D_1=\hat\pp
(S^2\cB)$.

\vspace{.5cm} (2) We use Lemma \ref{4.4}, \eqref{eqn:4.2} and
\eqref{eq4.4}. Note that $\Sigma_R^0$ does not intersect with
$\Omega_R$ and $\Delta_R$. Hence $D_2^0$ is the quotient of
$q^*\pp \Upsilon^{-1}(0)\git O(2)$ which is a subset of
$q^*\pp[\cE xt^1_{\mathcal{J}^{[n]}_0}(\cL_1,\cL_2)\oplus \cE
xt^1_{\mathcal{J}^{[n]}_0}(\cL_2,\cL_1)]\git O(2)$, by the action
of $PGL(N)$. The above are bundles over the restriction of
$$\pp\mathrm{Isom}(\cc^N,p_*\cL_1(k)\oplus p_*\cL_2(k))\git O(2)$$
to   $\mathcal{J}^{[n]}_0$. As in the proof of (1), observe that
$D_2^0$ is in fact the quotient of $q^*\pp\Upsilon^{-1}(0)$ by the
action of $PGL(N)\times O(2)$ since the actions commute. So we can
first take the quotient by the action of $PGL(N)$, then by the
action of $SO(2)$, and finally by the action of
$\zz_2=O(2)/SO(2)$. Since
$\pp\mathrm{Isom}(\cc^N,p_*\cL_1(k)\oplus p_*\cL_2(k))$ is a
principal $PGL(N)$-bundle, the quotient by $PGL(N)$ gives us
$$\pp\Upsilon^{-1}(0)\subset \pp[\cE
xt^1_{\mathcal{J}^{[n]}_0}(\cL_1,\cL_2)\oplus \cE
xt^1_{\mathcal{J}^{[n]}_0}(\cL_2,\cL_1)]$$ over
$\mathcal{J}^{[n]}_0$. The algebraic vector bundles $\cE
xt^1_{\mathcal{J}^{[n]}_0}(\cL_1,\cL_2)$ and $\cE
xt^1_{\mathcal{J}^{[n]}_0}(\cL_2,\cL_1)$ are certainly Zariski
locally trivial and in fact these bundles are dual to each other
by the Yoneda pairing $\Upsilon$ which is non-degenerate. In
particular, $\Upsilon^{-1}(0)$ is Zariski locally trivial.

Next we take the quotient by the action of $SO(2)\cong \cc^*$.
This action is trivial on the base $\mathcal{J}^{[n]}_0$ and
$SO(2)$ acts on the fibers. Hence $\pp\Upsilon^{-1}(0)/SO(2)$ is a
Zariski locally trivial subbundle of
$$\pp[\cE xt^1_{\mathcal{J}^{[n]}_0}(\cL_1,\cL_2)\oplus \cE
xt^1_{\mathcal{J}^{[n]}_0}(\cL_2,\cL_1)]\git \cc^*\cong \pp\cE
xt^1_{\mathcal{J}^{[n]}_0}(\cL_1,\cL_2)\times_{\mathcal{J}^{[n]}_0}
\pp\cE xt^1_{\mathcal{J}^{[n]}_0}(\cL_2,\cL_1)$$ over
$\mathcal{J}^{[n]}_0$ given by the incidence relations in terms of
the identification $$\pp\cE
xt^1_{\mathcal{J}^{[n]}_0}(\cL_1,\cL_2)\cong \pp\cE
xt^1_{\mathcal{J}^{[n]}_0}(\cL_2,\cL_1)^\vee.$$ Finally, $D_2^0$
is the $\zz_2$-quotient of $\pp\Upsilon^{-1}(0)/SO(2)$.

\vspace{.5cm}(3) By (a proof parallel to) \cite{og97} (1.7.10),
the intersection of $\Sigma_R^{ss}$ and $\Omega_R$ is smooth. By
Corollary \ref{4.3} and \eqref{4.1}, $\Delta_S$ is the blow-up of
$q^*\pp\mathrm{Hom}^{\omega}_2(W,T_{J^{[n]}\times \hat J})\git
SO(W)$ along $q^*\pp\mathrm{Hom}^{\omega}_1(W,T_{J^{[n]}\times
\hat J})\git SO(W)$. Hence $\Delta_S\git PGL(N)$ is the quotient
of
$$Bl_{q^*\pp\mathrm{Hom}^{\omega}_1(W,T_{J^{[n]}\times \hat J})}q^*\pp\mathrm{Hom}^{\omega}_2(W,T_{J^{[n]}\times \hat J})$$
by the action of $SO(W)\times PGL(N)$. By taking the quotient by
the action of $PGL(N)$ we get
$$Bl_{\pp\mathrm{Hom}^{\omega}_1(W,T_{J^{[n]}\times \hat J})}\pp\mathrm{Hom}^{\omega}_2(W,T_{J^{[n]}\times \hat J})$$
since $q$ is a principal $PGL(N)$-bundle. Next we take the
quotient by the action of $SO(W)$. Let
$\mathrm{Gr}^\omega(2,T_{J^{[n]}\times \hat J})$ be the relative
Grassmannian of isotropic 2-dimensional subspaces in
$T_{J^{[n]}\times \hat J}$ and let $\mathcal A$ be the
tautological rank 2 bundle on
$\mathrm{Gr}^\omega(2,T_{J^{[n]}\times \hat J})$. We claim
\begin{equation}\label{eq4.5}
Bl_{\pp\mathrm{Hom}^{\omega}_1(W,T_{J^{[n]}\times \hat
J})}\pp\mathrm{Hom}^{\omega}_2(W,T_{J^{[n]}\times \hat J})\git
SO(W)\cong \pp (S^2\mathcal A)
\end{equation}
which is a $\pp^2$-bundle over a $\mathrm{Gr}^\omega(2,2n)$-bundle
over $J^{[n]}$. It is obvious that the bundles are Zariski locally
trivial.

There are forgetful maps
\begin{eqnarray*}f:\pp\Hom(W,\cA)\to\pp\Hom^\omega_2(W,T_{J^{[n]}\times \hat J})\ \mbox{\rm
and}\\
f_1:\pp\Hom_1(W,\cA)\to\pp\Hom_1^\omega(W,T_{J^{[n]}\times \hat
J})\end{eqnarray*} where the subscript 1 denotes the locus of rank
$\le 1$ homomorphisms. Because the ideal of
$\pp\Hom_1^\omega(W,T_{J^{[n]}\times \hat J})$ pulls back to the
ideal of $\pp\Hom_1(W,\cA)$, $f$ lifts to
$$\hat f:Bl_{\pp \mathrm{Hom}_1(W,\cA)}\pp\Hom(W,\cA)\to Bl_{\pp\mathrm{Hom}^{\omega}_1(W,T_{J^{[n]}\times \hat J})}\pp\Hom^\omega_2(W,T_{J^{[n]}\times \hat J}).$$
This map is bijective (\cite[(3.5.1)]{og97}) and hence $\hat f$ is
an isomorphism. Now observe that the quotient $Bl_{\pp
\mathrm{Hom}_1(W,\cA)}\pp\Hom(W,\cA)\git SO(W)$ is isomorphic to
$\pp (S^2\mathcal A)$ where the quotient map is given by
$\alpha\mapsto \alpha\circ\alpha^t$. So we proved that
\begin{equation}\label{eq4.6}\Delta_S\git PGL(N)\cong \pp (S^2\mathcal
A).\end{equation} Finally $S\git PGL(N)$ is singular only along
$\Delta_S\git PGL(N)$ and the singularities are $\cc^{2n-1}/\{\pm
1\}$ by Luna's slice theorem \cite[(1.2.1)]{og97}. Since $D_3$ is
the exceptional divisor of the blow-up of $S\git PGL(N)$ along
$\Delta_S\git PGL(N)$, we conclude that $D_3$ is a
$\pp^{2n-2}$-bundle over $\pp (S^2\mathcal A)$.

\vspace{.5cm} (4) By Corollary \ref{4.3}, $\Sigma_S^s\cap\Omega_S$
is the exceptional divisor of the blow-up
$Bl_{q^*\pp\mathrm{Hom}^{\omega}_1(W,T_{J^{[n]}\times \hat
J})}q^*\pp\mathrm{Hom}^{\omega}(W,T_{J^{[n]}\times \hat J})\git
SO(W) $ and $\Sigma_T\cap \Omega_T$ is now the blow-up of the
exceptional divisor along the proper transform of
$$q^*\pp\mathrm{Hom}^{\omega}_2(W,T_{J^{[n]}\times \hat J})\git
SO(W).$$ Using the isomorphism \eqref{eq4.-2}, this is the
exceptional divisor of
$$q^*Bl_{\pp (S^2_1\cB)}\pp (S^2\cB)\to q^*\pp (S^2\cB)$$ over
$\mathrm{Gr}^\omega(3,T_{J^{[n]}\times \hat J})$. Since $q$ is a
principal $PGL(N)$-bundle, $D_1\cap D_2=\Sigma_T\cap \Omega_T\git
PGL(N)$ is the exceptional divisor of the blow-up $Bl_{\pp
(S^2_1\cB)}\pp (S^2\cB)$. As $\pp (S^2_1\cB)$ is a $\pp^2$-bundle
over $\mathrm{Gr}^\omega(3,T_{J^{[n]}\times \hat J})$, the
exceptional divisor is a $\pp^2\times \pp^2$-bundle over
$\mathrm{Gr}^\omega(3,T_{J^{[n]}\times \hat J})$. This is
obviously Zariski locally trivial.

\vspace{.5cm} (5) From the above proof of (3) it follows
immediately that $\Sigma_S^s\cap \Delta_S\git SO(W)$ is $\pp
(S^2_1\cA)$ and $D_2\cap D_3$ is a $\pp^{2n-2}$ bundle over $\pp
(S^2_1\cA)$.

\vspace{.5cm} (6) As in the above proof of (4), we start with
\eqref{eq4.1} and use the isomorphism \eqref{eq4.-2} to see that
$D_1\cap D_3$ is the proper transform of $\pp (S^2_2\cB)$ in the
blow-up $Bl_{\pp (S^2_1\cB)}\pp (S^2\cB)$. This is a Zariski
locally trivial $\pp^2\times\pp^2$-bundle over
$\mathrm{Gr}^\omega(3,T_{J^{[n]}\times \hat J})$.

\vspace{.5cm}(7) The description of $D_1\cap D_2\cap D_3$ follows
immediately from the proof of (4) and (6).

\vspace{.5cm}
 From the above descriptions, it is clear that $D_i$
($i=1,2,3$) are normal crossing smooth divisors. \qed \vspace{1cm}

 In order
to compute the stringy E-function of $\bM$ by using Kirwan's
desingularization $\widehat\bM$ and Definition \ref{def2.1}, we
also need the discrepancy divisor $K_{\widehat
\bM}-\rho^*K_{\bM}$.
\begin{proposition}\label{discrepancy divisor}
The discrepancy divisor of $\rho: {\widehat \bM}\to \bM$ is
\[(6n-1)D_1+(2n-2)D_2+(4n-2)D_3
\]
\end{proposition}
\begin{proof}
The proof is identical to that of \cite[(3.4.1)]{og97} and so we
omit the details.
\end{proof}
In particular, the singularities of $\bM=\bM_{2n}$ are terminal
for $n\ge 2$.
\begin{remark} Another way to prove Proposition \ref{discrepancy
divisor} is as follows. First observe as in \cite{og97} that
$\widehat{\bM}$ can be blown-down twice: $$\widehat{\bM}\to
\overline{\bM}\to \widetilde{\bM}$$ The first map is the
contraction of $D_3$ along the $\pp^2$-fiber and the second map is
the contraction of $D_1$ along the $\pp^5$-fiber (after the first
contraction $\hat \pp^5$ becomes $\pp^5$). The result of the two
contractions is also a desingularization $\nu:\widetilde{\bM}\to
\bM$. Since the singularities along $\Sigma$ are toric, it is easy
to compute the discrepancy along $D_2^0$ of $\nu$ which is
precisely $2n-2$ by toric geometry. It is not hard to check that
the pull-back of the closure of $D_2^0$ in $\widetilde{\bM}$ to
$\widehat{\bM}$ is $3D_1+D_2+2D_3$. From the well-known formula
\cite[II Ex. 8.5]{Hart}, we deduce that the discrepancy divisor
for $\rho$ is
$$(2n-2)(3D_1+D_2+2D_3)+5D_1+2D_3=(6n-1)D_1+(2n-2)D_2+(4n-2)D_3.$$\end{remark}


\section{Nonexistence of a crepant resolution}
\label{sec: Nonexistence of Crepant Resolution}

In this section we first find an expression for the stringy
E-function of the moduli space $\bM=\bM_{2n}$ with $n\geq 2$ by
using the detailed analysis of Kirwan's desingularization in
\S\ref{sec: Kirwan's desingularization}. Then we show that it
cannot be a polynomial, which proves Proposition \ref{prop:stringy
E-function test}.

By (\ref{eqn:stringy E-function}) and Proposition \ref{discrepancy
divisor}, the stringy E-function of ${\bf M}$ is given by
\begin{equation} \label{eqn:stringy E-function of M_c}
\begin{array}{lll}
& E({\bf M}^s;u,v)+ E(D_1^0;u,v)
{\textstyle\frac{1-uv}{1-(uv)^{6n}}
+E(D_2^0;u,v){\textstyle\frac{1-uv}{1-(uv)^{2n-1}}}} & \\
& +E(D_3^0;u,v) {\textstyle\frac{1-uv}{1-(uv)^{4n-1}}}
+E(D_{12}^0;u,v){\textstyle\frac{1-uv}{1-(uv)^{6n}}
\frac{1-uv}{1-(uv)^{2n-1}}} & \\
& +E(D_{23}^0;u,v){\textstyle\frac{1-uv}{1-(uv)^{2n-1}}
\frac{1-uv}{1-(uv)^{4n-1}}} & \\ &
+E(D_{13}^0;u,v){\textstyle\frac{1-uv}{1-(uv)^{4n-1}}
\frac{1-uv}{1-(uv)^{6n}}} &  \\
& +E(D_{123}^0;u,v){\textstyle\frac{1-uv}{1-(uv)^{6n}}
\frac{1-uv}{1-(uv)^{2n-1}}\frac{1-uv}{1-(uv)^{4n-1}}} &
.\end{array}
\end{equation}

We need to compute the Hodge-Deligne polynomials of $D^0_J$ for
$J\subset \{1,2,3\}$. Recall that for  a symplectic vector space
$(\cc^{2n},\omega)$, $\Gr^{\omega}(k,2n)$ denotes the Grassmannian
of $k$-dimensional subspaces of $\cc^{2n}$, isotropic with respect
to the symplectic form $\omega$ (i.e. the restriction of $\omega$
to the subspace is zero).

\begin{lemma}\label{lem:Hodge poly of Gr}\cite[Lemma 3.1]{CK04} For $k\leq n$,
the Hodge-Deligne polynomial of
$\Gr^\omega(k,2n)$ is
\[\prod_{1\leq i\leq k} \frac{1-(uv)^{2n-2k+2i}}{1-(uv)^i}. \]
\end{lemma}

From Lemma \ref{lem:Hodge poly of Gr} and Proposition
\ref{prop:analysis on exc}, we have the following corollary by the
additive and multiplicative properties of the Hodge-Deligne
polynomial.

\begin{corollary}\label{eqn:computation of stringy E-function}
$$ E(D_1;u,v) = \Bigl({\textstyle
\frac{1-(uv)^6}{1-uv}-\!\frac{1-(uv)^3}{1-uv}+\!\bigl(\frac{1-(uv)^3}{1-uv}\bigr)^2}\Bigr)
\! \times\!\!\! \prod_{1\leq i\leq 3}\! \Bigl({\textstyle
\frac{1-(uv)^{2n-4+2i}}{1-(uv)^i}}\Bigr)\!\! \times\! E(
J^{[n]}\!\!\times\!\! \hat J ;u,v),$$

$$E(D_3;u,v)   = {\textstyle
\frac{1-(uv)^{2n-1}}{1-uv}\cdot\frac{1-(uv)^3}{1-uv}} \times
\prod_{1\leq i\leq 2}\Bigl({\textstyle \frac{1-(uv)^{2n-2+2i}}
{1-(uv)^i}}\Bigr)\times E(J^{[n]}\times \hat J;u,v), $$

$$ E(D_{12};u,v)   = \Bigl({\textstyle
\frac{1-(uv)^3}{1-uv}}\Bigr)^2\times \prod_{1\leq i\leq
3}\Bigl({\textstyle \frac{1-(uv)^{2n-4+2i}}{1-(uv)^i}}\Bigr)
\times E(J^{[n]}\times \hat J;u,v), $$

$$ E(D_{23};u,v)  = {\textstyle
\frac{1-(uv)^{2n-1}}{1-uv}\cdot\frac{1-(uv)^2}{1-uv}}
\times\prod_{1\leq i\leq 2}\Bigl({\textstyle
\frac{1-(uv)^{2n-2+2i}}{1-(uv)^i}}\Bigr)\times E(J^{[n]}\times
\hat J;u,v), $$

$$E(D_{13};u,v)   = {\textstyle \frac{
1-(uv)^3}{1-uv}\cdot\frac{1-(uv)^{2n-2}}{1-uv}} \times
\prod_{1\leq i\leq 2}\Bigl({\textstyle
\frac{1-(uv)^{2n-2+2i}}{1-(uv)^i}}\Bigr)\times E(J^{[n]}\times
\hat J;u,v), $$

$$ E(D_{123};u,v)   ={\textstyle
\frac{1-(uv)^2}{1-uv}\cdot\frac{1-(uv)^{2n-2}}{1-uv}}\times
\prod_{1\leq i\leq 2}\Bigl({\textstyle
\frac{1-(uv)^{2n-2+2i}}{1-(uv)^i}}\Bigr)\times E(J^{[n]}\times
\hat J;u,v).$$
\end{corollary}
\begin{proof} The only thing that doesn't follow from the
multiplicative property of Hodge-Deligne polynomial is the
equations for $D_3$ and $D_{23}$ but this is a direct consequence
of the Leray-Hirsch theorem \cite[p.195]{V02I}.\end{proof}

For the E-polynomial of $D_2^0$ we have the following lemma.

\begin{lemma}\label{lem: Hodge Deligne poly of D02}
$ E(D_2^0;z,z) $ is divisible by $\frac{1-(z^2)^{2n-1}}{1-z^2}$.
\end{lemma}

\proof Note that $$I_{2n-1}=\{((x_i),(y_j))\in \pp^{2n-1}\times
\pp^{2n-1}\,|\, \sum_{i=0}^{2n-1} x_iy_i=0\} $$ and that it admits
a $\zz_2$-action interchanging $x_i$ and $y_i$. It is elementary
(\cite{GH78} p. 606) to see that
$$H^*(I_{2n-1};\qq)\cong \qq[a,b]/\langle a^{2n}, b^{2n},
a^{2n-1}+a^{2n-2}b+a^{2n-3}b^2+\cdots+b^{2n-1} \rangle$$ where $a$
(resp. $b$) is the pull-back of the first Chern class of the
tautological line bundle of the first (resp. second) $\pp^{2n-1}$.
The $\zz_2$-action interchanges $a$ and $b$. Let
$H^*(I_{2n-1};\qq)^{\pm}$ be the $\pm 1$-eigenspace of the
$\zz_2$-action in $H^*(I_{2n-1};\qq)$. The invariant subspace
$H^*(I_{2n-1};\qq)^{+}$ of $H^*(I_{2n-1};\qq)$ is generated by
classes of the form $a^ib^j+a^jb^i$. As a vector space
$H^*(I_{2n-1};\qq)$ is
\begin{equation}\label{eq3.0}\qq\text{-span}\{a^ib^j\,|\, 0\le i\le 2n-1, 0\le j\le
2n-2\}\end{equation} while the invariant subspace is
$$\qq\text{-span}\{a^ib^j+a^jb^i\,|\, 0\le i\le j\le 2n-2\}.$$
The index set $\{(i,j)\,|\, 0\le i\le j\le 2n-2\}$ is mapped to
its complement in $\{(i,j)\,|\, 0\le i\le 2n-1, 0\le j\le 2n-2\}$
by the map $(i,j)\mapsto (j+1,i)$. This immediately implies that
the Poincar\'e polynomial satisfies
\begin{equation}\label{eq3.3}
 P(I_{2n-1};z)=(1+z^2)P^+(I_{2n-1};z)
 \end{equation}
 where $P^{\pm}(I_{2n-1};z)=\sum (-1)^rz^r\dim
 H^r(I_{2n-1})^{\pm}$.
By \eqref{eq3.0}, we have
\begin{equation*}
 P(I_{2n-1};z)=\frac{1-(z^2)^{2n}}{1-z^2}\cdot
 \frac{1-(z^2)^{2n-1}}{1-z^2}.
 \end{equation*}
Because $1+z^2$ divides $\frac{1-(z^2)^{2n}}{1-z^2}$,
$\frac{1-(z^2)^{2n-1}}{1-z^2}$ also divides $P^+(I_{2n-1};z)$. By
\eqref{eq3.3}, $P^-(I_{2n-1};z)=z^2P^+(I_{2n-1};z)$ and hence
$\frac{1-(z^2)^{2n-1}}{1-z^2}$ also  divides $P^-(I_{2n-1};z)$.

Let
$$\psi:\mathcal{D}:=\pp\Upsilon^{-1}(0)/SO(2)\to
\mathcal{J}^{[n]}_0=(J^{[n]}\times \hat{J})\times
(J^{[n]}\times\hat{J})-\mathbf{\Delta}$$ be the Zariski locally
trivial $I_{2n-1}$-bundle in the proof of Proposition
\ref{prop:analysis on exc} (2). Recall that
$D_2^0=\mathcal{D}/\zz_2$. We have seen in the proof of
Proposition \ref{prop:analysis on exc} (2) that there is a
$\zz_2$-equivariant embedding
$$\imath:\mathcal{D}\hookrightarrow
\pp\cE
xt^1_{\mathcal{J}^{[n]}_0}(\cL_1,\cL_2)\times_{\mathcal{J}^{[n]}_0}
\pp\cE xt^1_{\mathcal{J}^{[n]}_0}(\cL_2,\cL_1)$$ where the
$\zz_2$-action interchanges $\cL_1$ and $\cL_2$.

Let $\lambda$ (resp. $\eta$) be the pull-back to $\mathcal{D}$ of
the first Chern class of the tautological line bundle over $\pp\cE
xt^1_{\mathcal{J}^{[n]}_0}(\cL_1,\cL_2)$ (resp. $\pp\cE
xt^1_{\mathcal{J}^{[n]}_0}(\cL_2,\cL_1)$). By definition,
$\lambda$ and $\eta$ restrict to $a$ and $b$ respectively. The
$\zz_2$-action interchanges $\lambda$ and $\eta$. By the
Leray-Hirsch theorem (\cite{V02I} p.195), we have an isomorphism
\begin{equation}\label{eq3.5} H^*_c(\mathcal{D})\ \
\cong \ H^*_c(\mathcal{J}^{[n]}_0 )\otimes
H^*(I_{2n-1}).\end{equation} As the pull-back and the cup product
preserve mixed Hodge structure, \eqref{eq3.5} determines the mixed
Hodge structure of $H^*_c(\mathcal{D})$. The $\zz_2$-invariant
part is
\begin{equation*}\label{eq3.6}
H^*_c(\mathcal{D})^+\cong \left( H^*_c( \mathcal{J}^{[n]}_0
)^+\otimes H^*(I_{2n-1})^+\right) \oplus \left(
H^*_c(\mathcal{J}^{[n]}_0 )^-\otimes H^*(I_{2n-1})^-\right)
\end{equation*} where the superscript $\pm$ denotes the $\pm
1$-eigenspace of the $\zz_2$-action. Because $H^*_c(D_2^0)\cong
H^*_c(\mathcal{D}/\zz_2)\cong H^*_c(\mathcal{D})^+$ (\cite{Gr57}
Theorem 5.3.1 and Proposition 5.2.3), $E(D_2^0;u,v)$ is equal to
\begin{equation*}\label{eq4.7}
 E^+(\mathcal{D};u,v)=E^+( \mathcal{J}^{[n]}_0;u,v)E^+(I_{2n-1};u,v)
 +E^-(\mathcal{J}^{[n]}_0;u,v)E^-(I_{2n-1};u,v).
 \end{equation*}
where $E^\pm(Y;u,v)=\displaystyle\sum_{p,q}\sum_{k\geq0} (-1)^k
h^{p,q}(H^k_c(Y)^\pm) u^pv^q$. Since the smooth projective variety
$I_{2n-1}$ has pure Hodge structure,
$$E^+(I_{2n-1};z,z)=P^+(I_{2n-1};z)\ \mbox{and}\
E^-(I_{2n-1};z,z)=P^-(I_{2n-1};z) .$$ As
$\frac{1-(z^2)^{2n-1}}{1-z^2}$ divides $P^\pm(I_{2n-1};z)$, it
divides $E(D_2^0;u,v)$ as well. \qed\\

\textit{Proof of Proposition \ref{prop:stringy E-function test}.}

Let us prove that \eqref{eqn:stringy E-function of M_c} cannot be
a polynomial. Let
$$S(z)=E_{st}({\bf M};z,z)-E({\bf M}^s;z,z)-\frac{1-z^2}{1-(z^2)^{2n-1}}E(D_2^0
;z,z).$$ It suffices to show that $S(z)$ is not a polynomial for
all $n\geq2$ because $E({\bf M}^s;z,z)$ and
$\frac{1-z^2}{1-(z^2)^{2n-1}}E(D_2^0 ;z,z)$ are  polynomials by
Lemma \ref{lem: Hodge Deligne poly of D02}.

Express the rational function $S(z)$ as
$$\frac{N(z)}{(1-(z^2)^{2n-1})(1-(z^2)^{4n-1})(1-(z^2)^{6n})}.$$
By direct computation using (\ref{eqn:stringy E-function of M_c})
and Corollary \ref{eqn:computation of stringy E-function}, $N(z)$
modulo $1-(z^2)^{2n-1}$ is congruent to
\begin{equation}\label{eqn:denumerator modulo}
\begin{array}{lll}(1-z^2)^2(1-(z^2)^{4n-1})\times\Bigl({\textstyle
\frac{1-(z^2)^3}{1-z^2}}\Bigr)^2\times \prod_{1\leq i\leq
3}\Bigl({\textstyle \frac{1-(z^2)^{2n-4+2i}}{1-(z^2)^i}}\Bigr)
\times P(\mathfrak J;z)&& \\ -(1-z^2)^2(1-(z^2)^{4n-1})\times
{\textstyle
\frac{1-(z^2)^2}{1-z^2}\cdot\frac{1-(z^2)^{2n-2}}{1-z^2}}\times
\prod_{1\leq i\leq 2}\Bigl({\textstyle
\frac{1-(z^2)^{2n-2+2i}}{1-(z^2)^i}}\Bigr) \times P(\mathfrak J;z)
&& \\ -(1-z^2)^2(1-(z^2)^{6n})\times {\textstyle
\frac{1-(z^2)^2}{1-z^2}\cdot\frac{1-(z^2)^{2n-2}}{1-z^2}}\times
\prod_{1\leq i\leq 2}\Bigl({\textstyle
\frac{1-(z^2)^{2n-2+2i}}{1-(z^2)^i}}\Bigr) \times P(\mathfrak J;z)
&& \\ + (1-z^2)^3\times{\textstyle \frac{1-(z^2)^2}{1-z^2}
\cdot\frac{1-(z^2)^{2n-2}}{1-z^2}}\times \prod_{1\leq i\leq
2}\Bigl({\textstyle
\frac{1-(z^2)^{2n-2+2i}}{1-(z^2)^i}}\Bigr)\times P(\mathfrak J;z)
&& \end{array} \end{equation} where $\mathfrak
J:=J^{[n]}\times\hat J$.

All we need to show is that the numerator $N(z)$ is not divisible
by the denominator
$(1-(z^2)^{2n-1})(1-(z^2)^{4n-1})(1-(z^2)^{6n})$. We write
(\ref{eqn:denumerator modulo}) as a product $s(t)\cdot P(\mathfrak
J;z)$ for some polynomial $s(t)$ with $t=z^2$. For the proof of
Proposition \ref{prop:stringy E-function test} for $n\geq 3$ (the
$n=2$ case will be treated separately), it suffices to prove the
following:
\begin{enumerate} \item if $n+1$ is not divisible by 3, then $1-z^2$
is the GCD of $1-(z^2)^{2n-1}$ and $s(z^2)$, and
$\frac{1-(z^2)^{2n-1}}{1-z^2}$ does not divide $P(\mathfrak J;z)$;
\item if $n+1$ is divisible by 3, then $1-(z^2)^3$ is the GCD of
$1-(z^2)^{2n-1}$ and $s(z^2)$, and
$\frac{1-(z^2)^{2n-1}}{1-(z^2)^3}$ does not divide $P(\mathfrak
J;z)$.
\end{enumerate}

For (1), suppose $n+1$ is not divisible by 3. From
(\ref{eqn:denumerator modulo}), $s(t)$ is divisible by $1-t$. We
claim that $s(t)$ is not divisible by any irreducible factor of
$\frac{1-t^{2n-1}}{1-t}$, i.e. for any root $\alpha$ of
$1-t^{2n-1}$ which is not 1, $s(\alpha)\neq 0$. Using the relation
$\alpha^{2n-1}=1$, we compute directly that
\begin{equation}\label{eqn:bar s} s(\alpha)={\textstyle
-\frac{\alpha(1-\alpha^{-1}){(1-\alpha^3)}^2}{1+\alpha}},
\end{equation} which is not 0 because 3 does not divide $2n-1$.

Next we check that $\frac{1-(z^2)^{2n-1}}{1-z^2}$ does not divide
$P(\mathfrak J;z)$. Note that
\begin{equation}\label{eqnfinal002}P(\mathfrak
J;z)=P(J^{[n]};z)P(\hat{J};z)=(1-z)^4P(J^{[n]};z)\end{equation}
and hence it suffices to show that $\frac{1-(z^2)^{2n-1}}{1-z^2}$
does not divide $P(J^{[n]};z)$. We put $P(J^{[n]};z)= \sum_{0\leq
i\leq 4n} a_i z^i$ and write
\begin{eqnarray*}
\lefteqn{ \sum_{0\leq i\leq 4n} a_i z^i  = a_0 +a_1z
+a_2z^2+...+a_{4n-4}z^{4n-4}} &&   \\ && +a_{4n-3}z
\Bigl(\frac{1-(z^2)^{2n-1}} {1-z^2}-\sum^{2n-3}_{i=0}(z^2)^i
\Bigr) + a_{4n-2}(z^{4n-2}-1)+a_{4n-2}\nonumber
\\
& & + a_{4n-1}z(z^{4n-2}-1)+a_{4n-1}z+ a_{4n
}z^2(z^{4n-2}-1)+a_{4n}z^2.\nonumber
\end{eqnarray*}
We see from this that $P(J^{[n]};z)$ is divisible by
$\frac{1-(z^2)^{2n-1}}{1-z^2}$ if and only if
\begin{equation}\label{eqfinal000}\begin{array}{lll}
\lefteqn{ a_0 +a_1z +a_2z^2+...+a_{4n-4}z^{4n-4}} &&   \\ &&
+a_{4n-3}z \Bigl(-\sum^{2n-3}_{i=0}(z^2)^i \Bigr) + a_{4n-2}+
a_{4n-1}z+a_{4n}z^2\end{array}
\end{equation}
is divisible by $\frac{1-(z^2)^{2n-1}}{1-z^2}$. Since
\eqref{eqfinal000} is of degree $\le 4n-4$, it is divisible by
$\frac{1-(z^2)^{2n-1}}{1-z^2}$ if and only if \eqref{eqfinal000}
is a constant multiple of $\frac{1-(z^2)^{2n-1}}{1-z^2}$. If this
were true then the coefficient of $z$ must be zero, i.e.
$a_1-a_{4n-3}+a_{4n-1}=0$. By the Poincar\'{e} duality
$a_1-a_{4n-3}+a_{4n-1}= a_1-a_3 +a_1 $. This value is not 0
because $a_1=-b_1(J^{[n]})=-4$ and $a_3=-b_3(J^{[n]})=-40$ for
$n\geq3$ by G\"{o}ttsche's formula \cite{Go90}:
\begin{equation}\label{eqn:Betti for X[n]} \sum_{n\geq
0}P(J^{[n]};z)t^n=\prod_{k\geq1}
\prod_{i=0}^{4}(1-z^{2k-2+i}t^k)^{(-1)^{i+1}b_i(J)}.
\end{equation}

For (2), suppose 3 divides $n+1$ and $n\neq 2$. Then from
(\ref{eqn:bar s}), $(1-t^3)$ divides $s(t)$. More precisely, for a
third root of unity $\alpha$, $s(\alpha)=0$. On the other hand, if
$\alpha$ is a root of $1-t^{2n-1}$ but not a third root of unity
then we can observe that $s(\alpha)\neq 0$ by (\ref{eqn:bar s}).
Since every root of $1-t^{2n-1}$ is a simple root, any irreducible
factor of $\frac{1-t^{2n-1}}{1-t^3}$ does not divide $s(t)$.

We next check that the polynomial
$\frac{1-(z^2)^{2n-1}}{1-(z^2)^3}$ does not divide $P(\mathfrak
J;z)$. Again by \eqref{eqnfinal002}, it suffices to show that
$\frac{1-(z^2)^{2n-1}}{1-(z^2)^3}$ does not divide $P(J^{[n]};z)$.
Let $P(J^{[n]};z)=\sum_{0\leq i\leq 4n} a_i z^i$ and write
$z^{4n-8}$ as
$D(z)=\frac{1-(z^2)^{2n-1}}{1-(z^2)^3}-\sum_{i=0}^{\frac{2n-7}3}(z^2)^{3i}$.
Then we have
\begin{eqnarray*}\label{eqn:p(t) when 3 divides n}
\lefteqn{ \sum_{0\leq i\leq 4n} a_i
z^i=a_0+a_1z+...+a_{4n-8}z^{4n-8} } &&
\\ && + a_{4n-7}zD(z)+a_{4n-6}z^2D(z)+...+a_{4n-3}z^5D(z)  \nonumber
\\ && +
a_{4n-2}(z^{4n-2}-1)+a_{4n-2}+a_{4n-1}z(z^{4n-2}-1)+a_{4n-1}z \nonumber \\
&& + a_{4n}z^2(z^{4n-2}-1)+a_{4n}z^2. \nonumber \end{eqnarray*}
Therefore, $P(J^{[n]};z)$ is divisible by
$\frac{1-(z^2)^{2n-1}}{1-(z^2)^3}$ only if
\begin{equation}\label{eqnfinal001}\begin{array}{lll}
\lefteqn{ a_0+a_1z+...+a_{4n-8}z^{4n-8} } &&
\\ && + (a_{4n-7}z+a_{4n-6}z^2+...+a_{4n-3}z^5)(-\sum_{i=0}^{\frac{2n-7}3}(z^2)^{3i})
\\ &&
+a_{4n-2}+a_{4n-1}z +a_{4n}z^2\end{array}\end{equation} is
divisible by $\frac{1-(z^2)^{2n-1}}{1-(z^2)^3}$. Since
\eqref{eqnfinal001} is of degree $\le 4n-8$, it is divisible by
$\frac{1-(z^2)^{2n-1}}{1-(z^2)^3}$ if and only if it is a constant
multiple of $\frac{1-(z^2)^{2n-1}}{1-(z^2)^3}$. If this were true
the coefficient of $z$ must be zero, i.e.
$a_1-a_{4n-7}+a_{4n-1}=0$. By the Poincar\'{e} duality
$a_1-a_{4n-7}+a_{4n-1}= a_1-a_7 +a_1 $. This value is not zero
because $a_1 =-4$ and $a_7=-b_7(J^{[n]})\leq -196$ for $n\geq3$ by
direct computation using G\"{o}ttsche's formula again.

The case of $n=2$ remains to be proved. We show that $N(z)$ is not
divisible by $1-(z^2)^{6n}=1-(z^2)^{12}$. By direct computation
using (\ref{eqn:stringy E-function of M_c}) and Corollary
\ref{eqn:computation of stringy E-function}, we have
\begin{equation*}\begin{array}{lll}\lefteqn{ N(z)= (1-z^2 )
(1+ z^2 ) (1+(z^2)^2) ((z^2)^2- z^2 +1)  (1+ z^2 +(z^2)^2)} && \\
&&  \times  \bigl((z^2)^{12}+3 (z^2)^{11}+3 (z^2)^{10}+2 (z^2)^9+2
(z^2)^8+3 (z^2)^7+3
(z^2)^6\\
&& +(z^2)^5+(z^2)^3+(z^2)^2+1\bigr)\times P(\mathfrak
J;z).\end{array}\end{equation*} By G\"{o}ttsche's formula, we also
have
$$P(\mathfrak J;z)= (1-4z+13 z^2-32 z^3+44 z^4-32 z^5+13 z^6-4
z^7+z^8 )(1-z)^4.$$ By plugging in a primitive root of
$z^{24}-1=0$, it is easy to check that $N(z)$ is not divisible by
$1-(z^2)^{12}$.

Therefore, $E_{st}(\bM;z,z)$ is not a polynomial for any $n\ge 2$.
 \qed \\

\begin{remark} The sum of second Chern class together with the determinant map
give us a morphism $$a:\bM\to J\times \mathrm{Pic}^0(J).$$ Let
$\cM=a^{-1}(0,0)$. Like $\bM$, $\cM$ is a singular projective
variety equipped with a holomorphic symplectic form on the smooth
part. One may ask if there is a crepant resolution of $\cM$. It is
easy to modify our proof to show that there is no crepant
resolution of $\cM$ (and therefore no symplectic
desingularization) either. We leave the details to the reader.
\end{remark}


\bibliographystyle{amsplain}

\end{document}